\documentstyle[11pt]{article}

\newcommand{\funbi}[5]
{
{#1} \left\{
\begin{array}{ll}
{#2} & \mbox{{#3}} \\
{#4} & \mbox{{#5}}
\end{array}
\right.
}

\newtheorem{theo}{Theorem}

\newtheorem{coro}{Corollary}

\newcommand{\proof}{\vspace{0.00cm}{\bf Proof: }}

\newcommand{\proofend}{{\bf Q.E.D.}\vspace{0.0cm}}

\title{
Collective Property of Numbers and Its Mathematical Refutation
}

\author{
\vspace{0.00cm}
Guang-Liang Li 
\footnote{The corresponding author.}
\hspace{2cm}
Victor O. K. Li
\\
Department of Electrical and Electronic Engineering\\
The University of Hong Kong\\
Rm. 601, Chow Yei Ching Bldg., 
Pokfulam Rd., Hong Kong \\
\{glli,vli\}@eee.hku.hk, Phone: (852)2857 8495,
Fax: (852)2559 8738
}

\begin{document}
\maketitle




\begin{abstract}
{\em 
A number has the ``collective'' property if
the number is the greatest lower bound of
a bounded, strictly decreasing sequence on the
real line. We prove that numbers
with the collective property constitute an empty set.
}
\end{abstract}

The properties 
of numbers for   
counting, 
calculating, or measuring 
are conceivable without
considering ``all numbers''. 
In mainstream mathematics, 
``all numbers'' or
``numbers as a whole'' is a notion formulated
based on the notion of ``set'', and
a number also has some ``collective'' 
properties different
from the properties for counting, calculating, and measuring.
A specific collective property of a number is necessary to
the definition of ``bounded, strictly decreasing sequence''.
To define ``bounded, strictly decreasing sequence'', one has to define 
a number, which is not a term of
the sequence, as the greatest lower bound 
of the sequence.
The property of the number as the greatest lower bound may not
be necessary when the number is used for
counting, calculating or measuring.
In the following, the collective property of a number
means ``the number is the greatest lower bound of a bounded, strictly
decreasing sequence on the real line''.

Due to the 
failure to
prove the consistency of any
formal system involving numbers
with the collective property \cite{Ebbi},
and also due to the fact that the definition of
the greatest lower bound is circular and hence is not
on a logically secure ground \cite{Eppl,Kush}, 
in mainstream mathematics, one has to accept the collective property
based on a philosophical belief (realism or Platonism). 
As a sharp contrast,  
constructive mathematicians
reject the collective property \cite{Bees,Bish,Bish2,Brid,Kush}.
However, the rejection is also
based on a philosophical belief
(constructivism or intuitivism), 
perhaps for lack of a convincing way
to argue against the collective property mathematically.
By proving the following theorem, 
we present a mathematical refutation of the collective property. 
In other words, we prove that numbers
with the collective property constitute an empty set and
hence are devoid of any meaning in mathematics.

\begin{theo}  
\label{th-1}
Denote by $E$ a bounded set 
consisting of all terms of an arbitrarily given,
strictly decreasing sequence 
(with a greatest lower bound 
$a \not \in E$)
on the real line.
The above definition of $E$ implies a contradiction.
\end{theo}

\proof 
For each point $x \in (-\infty, \infty)$, define

\begin{equation}
\label{eq-f}
\funbi 
{f(x) =}{1,}{$x \in E$;} 
{0,}{otherwise.}
\end{equation}

By (\ref{eq-f}), $f(x) = 0$ does not hold for all $x$.
Let $P(x)$ represent ``$f(z) = 0$ for each $z \leq x$''. 
Let $v$ be a real-valued variable with the following
properties. (I) $v$ assumes values greater than or equal to 
$x$ with $P(x)$
in the order from left to right
as the points appear on the real line. 
(II) In the order mentioned in (I), 
the set of values
of $v$ may have a last element
defined by
a condition, and the value of $v$ can be
any point   
less than or equal to the
last value. (III) If $v$ does not have a last value, then
the value of $v$ can be any point on the real line. 
If $v$ runs from any $x_0$ with 
$P(x_0)$ towards the right, i.e.,
if $v$ increases from $x_0$,
then $v$ reaches eventually
a unique $x_1$ with $f(x_1) = 1$ and
$f(x) = 0$ for each $x < x_1$.
To see this, we define $v|x'\rightarrow x''$ to mean
``if $P(x')$ holds, then let $v$ take on 
each $x$ with $x \geq x'$ 
as its value, and
let $x''$ be the last value of $v$ if
$f(v) = 1$ where the value of $v$ is $x''$''. 
Note the difference between 
$f(v) = 1$ as
the condition to define the last value of $v$ 
and $f(x) = 1$, where $x$ is a point in $E$.

Assume $P(x')$ holds. Then
$v|x'\rightarrow x''$ determines a nonempty
subset of values of $v$. Denying ``$x''$ is the
last value of $v$'' amounts to 
negating ``$f(v) = 1$ where the value of $v$ is $x''$''
in the definition of $v|x'\rightarrow x''$. The negation
asserts ``$f(v) = 0$ where the value of $v$ can be
any $x \geq x'$''. 
This contradicts (\ref{eq-f}).
Since $x''$ is a value of $v$, no $x$ with 
$x < x''$ is the last value of $v$. So
$f(v) = 0$ where the value of $v$ can be any  
$x < x''$.
In other words, $x''$ is unique. 
The uniqueness of $x''$ is independent of $x'$ so long as
$P(x')$ holds.
By the definition of $v|x'\rightarrow x''$, 
given $P(x')$, if
$v|x'\rightarrow x''$, and if $v|x'\rightarrow y''$,
then $x'' = y''$.  
Consider two arbitrarily given points 
$x'$ and $y'$ with $x' \not = y',
P(x'),  P(y'),
v|x'\rightarrow x''$, and $v|y'\rightarrow y''$. 
If $x''\not = y''$, 
then without
loss of generality, assume $y' < x' < x'' < y''$, which
implies $v|y' \rightarrow x''$.
This contradicts  $v|y'\rightarrow y''$. So $x'' = y''$.
Consequently, for all $x$ with $P(x)$,
there is a unique $x_1$ with $v|x\rightarrow x_1$. 

The above
argument has a clear physical meaning.
Negating the argument then leads to an absurdity.
Let $r$ represent the position of
a point particle moving along the real line
towards the right from an initial position
$x_0$ with $P(x_0)$.
The particle moves at $r$ if and only if
$P(r)$ holds, i.e.,
the speed of the particle at $r$ is positive
if and only if $P(r)$ holds. Then the
particle stops at $x_1$ with
$f(x_1) = 1$ and never moves again.
If this is not the case, then the particle either
is motionless with a positive speed, or
moves at a speed equal to 0. 
This violates the basic physical law and hence is absurd.

If $x_1$ does not exist, then
$v$ runs towards $\infty$, which implies  
$f(x) = 0$ at each $x$.
This contradicts (\ref{eq-f}) as we have shown already.
However, if $x_1$ exists, then
$x_1 = \min E \in E$,
which contradicts the definition of $E$.
We see a contradiction in either case above. 
\proofend

Clearly, if we detach the collective property from numbers,
then any contradiction or absurdity caused by the
collective property disappears.
The collective property comes into being by the definition of the
greatest lower bound of $E$. The contradiction 
implied by the
definition of $E$ suggests that the definition of the
greatest lower bound is problematic. 
This is indeed the case.
The definition of the greatest lower bound is circular 
\cite{Eppl,Kush}, and
violates Russell's vicious circle principle:
No totality can contain members defined in terms of the
totality itself \cite{Guin}. Definitions
violating this principle are called non-predicative
(or impredicative). According to Russell, non-predicative
definitions
cause contradictions, such as the well-known 
Russell's paradox \cite{Eppl}. To avoid contradictions, 
Russell proposed the vicious circle principle.
Now let us take a look at 
the definition of the greatest lower bound.

A real number $y$ is a lower bound of an 
infinite set $S$ of real numbers, 
if $y \leq  x$ for each 
$x$ in $S$. Let $L$ be the set of all lower bounds of $S$. 
The greatest lower bound $l$ of $S$ 
is a member of $L$ with $l \geq  y$ 
for each $y$ in $L$. 
But $l$ may not necessarily be in $S$. 
Clearly, the definition of $l$ is circular and non-predicative: 
To define $l$, one has to define a totality $L$ first, 
in which $l$ itself is a member.
The contradiction implied by $E$ 
is directly due to the circular, 
non-predicative definition
of the greatest lower bound, which
attaches the
collective property to a point $a \not\in E$, and hence
excludes $\min E$ from $E$.
In fact, the definition of $E$ requires one to exclude $\min E$
from $E$. Achieved by
the circular, 
non-predicative definition
of the greatest lower bound,
such exclusion is logically invalid \cite{Eppl,Guin,Kush}.
As shown below,
the definition of the greatest lower bound of $E$ 
is merely a result of an artificial choice rather than a logical
necessity. The choice also causes an absurdity.

Denote by ${\overline N}(x)$ the number of different values
taken on by $f$ on $(-\infty, x]$, and $N'(x)$
the number of different values
taken on by $f$ on $(-\infty, x)$. 
Both ${\overline N}$ and $N'$ are
non-decreasing functions on $(-\infty, \infty)$.
If $f$ takes on only one value on $(-\infty, x]$,
then ${\overline N}(x) = N'(x) = 1$.
For any $x \in (-\infty, \infty)$, ${\overline N}(x) \geq N'(x)$.
Write 
$\Delta N(x) = {\overline N}(x) - N'(x)$.
So 

\begin{equation}
\label{eq-N}
{\overline N}(x) = N'(x) + \Delta N(x)
\end{equation}
with 

\begin{equation}
\label{eq-lim}
\lim_{x \rightarrow -\infty}{\overline N}(x) = 1. 
\end{equation}
Since $f$ has only two different values, 
$\Delta N(x) \in \{0, 1\}$. 

For each $x \in (-\infty, \infty)$,
all eligible (i.e., consistent with (\ref{eq-N}) and (\ref{eq-lim}))
combinations of 
$N'(x), \max\{f(y): y < x\}, \Delta N(x), 
{\overline N}(x)$, and $\max\{f(y): y \geq x\}$ give 
the conditions $A(x), B(x), C(x)$ and $D(x)$ below,
such that either
$(-\infty, \infty) = \{x: A(x)\}$
or $(-\infty, \infty) = \{x: B(x)\}\cup\{x: C(x)\}\cup\{x: D(x)\}$.
\begin{center}
\begin{tabular}{c|c|c|c|c|c} 
   & $N'(x)$ & $\max\{f(y): y < x\}$ & $\Delta N(x)$ & 
${\overline N}(x)$ & $\max\{f(y): y \geq x\}$\\ \hline
$A(x)$ &   1     &        0           &     0         &
       1     &         0\\ \hline

$B(x)$ &   1     &        0           &     0         &
       1     &            1 \\ \hline

$C(x)$ &   1     &        0           &     1         &
       2     &            1 \\ \hline

$D(x)$ &   2     &        1           &     0         &
       2     &            1 \\ \hline

\end{tabular}
\end{center}

If $\{x: A(x)\} = (-\infty, \infty)$,
then
$\{x: B(x)\} = \{x: C(x)\} = \{x: D(x)\} = \emptyset$. 
In other words, 0 is the
only value of $f$.
This is prohibited by (\ref{eq-f}). 
Negating ``$\{x: A(x)\} = (-\infty, \infty)$'',
we obtain
$\{x: B(x)\}\cup\{x: C(x)\}\cup\{x: D(x)\} = (-\infty, \infty)$,
where
none of $\{x: B(x)\}, \{x: C(x)\}$, and $\{x: D(x)\}$
is empty.
As can be readily seen,
$\{x: B(x)\} = (-\infty, x_1), \{x: C(x)\} = \{x_1\}$, and
$\{x: D(x)\} = (x_1, \infty)$ with $\Delta N(x_1) = 1$.

Consider two arbitrarily given points
$x_0 < x'$ with ${\overline N}(x') > {\overline N}(x_0)$.
We say ``${\overline N}$ increases'' if
there are points $y < x$ such that 
${\overline N}(y) < {\overline N}(x)$.
We say ``${\overline N}$ increases
at $x$'' if ${\overline N}(x) > {\overline N}(y)$ for any
$y < x$. 
Since ${\overline N}$ increases 
on $[x_0, x']$,
${\overline N}$      
increases  
at a unique $x_1 \in [x_0, x']$,
which implies $\Delta N(x_1) = 1$.
The increase in ${\overline N}$ can then be expressed by

\[
2 = {\overline N}(x') = {\overline N}(x_0) + \int_{-\infty}^{\infty} 
\delta(x - x_1)dx = 1 + 1 
\]
where $\delta$ is the Dirac delta function.

In mainstream mathematics, one
denies that the negation of 
``$\{x: A(x)\} = (-\infty, \infty)$'' is
``$\{x: B(x)\}\cup\{x: C(x)\}\cup\{x: D(x)\} = (-\infty, \infty)$
with nonempty $\{x: B(x)\}, \{x: C(x)\}$, and $\{x: D(x)\}$'' 
by choosing,
without any logically justified reason, to
negate ``$\Delta N(x_1) = 1$''.
The denial implies 
not only the definition of the greatest lower bound of
$E$ but also an absurdity: 
``${\overline N}$ does not increase at any point on
the real line and 
${\overline N}$ increases''.
Using the definition of the greatest lower bound of
$E$ to explain away
the absurdity is 
logically invalid.

Actually, any formal definition or 
informal notion of even $\{n: n = 1, 2, \cdots\}$
is circular and non-predicative \cite{Myhi,Pars}. 
A formal definition or informal notion of
$\{n: n = 1, 2, \cdots\}$ is the basis of mainstream
mathematics. Again, in mainstream mathematics,
the conception of $\{n: n = 1, 2, \cdots\}$ 
is not on a
logically secure ground, and the
only reason for one to accept 
$\{n: n = 1, 2, \cdots\}$ is   
a philosophical belief 
\cite{Myhi}.

Let $P_1(E)$ stand for ``$E$
is a bounded, strictly decreasing 
sequence on the real line'', 
$P_2(1/{\bf N})$ for ``$1/{\bf N}$ is $\{1/n: n = 1, 2, \cdots\}$
or any infinite subset of $\{1/n: n = 1, 2, \cdots\}$'', 
and $P_3({\bf N})$ for ``${\bf N}$ is $\{n: n = 1, 2, \cdots\}$
or any infinite subset of $\{n: n = 1, 2, \cdots\}$''. 
The following results are immediate from
Theorem \ref{th-1}.

\begin{coro}
\label{co-1}
(i) $\{E: P_1(E)\} = \emptyset$, which implies
(with $E = \{1/n: n = 1, 2, \cdots\}$)
$\{1/{\bf N}: P_2(1/{\bf N})\} = \emptyset$, and
(due to the correspondence between
$\{1/n: n = 1, 2, \cdots\}$ and $\{n: n = 1, 2, \cdots\}$)
$\{{\bf N}: P_3({\bf N})\} = \emptyset$.
(ii) $\{x: x\; \mbox{is a number with the collective property}\}
= \emptyset$.
\end{coro}

By Corollary \ref{co-1}, 
$\{n: n = 1, 2, \cdots\}$ is meaningless, 
though $\{1, 2, \cdots, n\}$ is meaningful
for any given positive integer $n$. 
Also by Corollary \ref{co-1}, no number has the
collective property. In other words, the
definition of the greatest lower bound 
defines a concept without any instance. In particular,
defining numbers based on
the collective property, 
such as the definition of
irrational numbers, actually defines nothing.
Consequently,
mathematical reasoning based on a number
system with the collective property may
yield concepts with an empty set of instances
and hence may be misleading. 
Nevertheless, any number with a
decimal expression of a finite length still possesses
the properties for counting, calculating and
measuring. 
As shown by
constructive mathematics, numbers without the
collective property
are sufficient for scientific and engineering
applications.

\end{document}